\documentclass[a4paper,11pt,reqno]{article}
\usepackage[T1]{fontenc}
\usepackage{courier}
\usepackage[scaled=0.92]{helvet}

\usepackage{amscd,amsmath,amssymb,amsthm,amsfonts,epsfig,graphics,accents}

\usepackage[english]{babel}
\usepackage{amsthm,amsfonts,amssymb,amsmath}
\usepackage{hyperref}
\hypersetup{
colorlinks=true,
citecolor=red,
linkbordercolor={1 1 1},
citebordercolor={1 1 1},
}

\usepackage{mathrsfs}

\usepackage{xcolor} 

\theoremstyle{plain}
\newtheorem{satz}{Theorem}[section]
\newtheorem{prop}[satz]{Proposition}
\newtheorem{cor}[satz]{Corollary}
\newtheorem{lem}[satz]{Lemma}

\theoremstyle{definition}

\newtheorem{rem}[satz]{Remark}

\newcommand{\mx}{\mbox}
\newcommand{\rw}{\rightarrow}

\newcommand{\ml}{\mathcal}

\newcommand{\pl}{\partial}

\newcommand{\x}{\times}

\newcommand{\beq}[1]{\begin{equation} \label{#1}}
\newcommand{\eeq}{\end{equation}}
\newcommand{\beqar}{\[ \begin{array}{rcl}}
\newcommand{\eeqar}{\end{array} \]}

\newcommand{\jj}{^{(j)}}
\newcommand{\ift}{^{(\infty)}}

\providecommand{\RR}{\mathbb{R}}
\providecommand{\CC}{\mathbb{C}}

\providecommand{\NN}{\mathbb{N}}
\newcommand{\snorm}[2]{\left| #1\right|_{#2}}
\newcommand{\norm}[2]{\left \lVert#1 \right\rVert_{#2}}

\newcommand{\lie}[1]{\mathcal{L}_{\chi^{(#1)}}}

\newcommand{\ff}[1]{(\ref{#1})}

\newcommand{\Heq}[2]{\overset{\left(#1\right)}{\underset{}{#2}}}

\usepackage[utf8]{inputenc}

\hypersetup{pdfauthor={Fortunati, Wiggins},%
            pdftitle={Integrability and strong normal forms...},%
            pdfsubject={},%
            pdfkeywords={Non-autonomous systems, normal forms, aperiodic time dependence.}%
}

\DeclareMathOperator{\diag}{diag}

\DeclareMathOperator{\id}{Id}

\makeatletter
\renewcommand*{\@fnsymbol}[1]{\ensuremath{\ifcase#1\or *\or (a)\or (b)\or \else \fi}}
\makeatother
\title{\LARGE{\textbf{Integrability and strong normal forms for non-autonomous systems in a neighbourhood of an equilibrium\footnote{This research was supported by ONR Grant No.~N00014-01-1-0769 and MINECO: ICMAT Severo Ochoa project SEV-2011-0087.}}}}
\date{}
\author{%
Alessandro Fortunati\thanks{E-mail: alessandro.fortunati@bristol.ac.uk} \and 
Stephen Wiggins\thanks{E-mail: s.wiggins@bristol.ac.uk}
\bigskip \\
School of Mathematics, University of Bristol, Bristol BS8 1TW, United Kingdom
}

\begin{document}

\maketitle

\begin{abstract}
The paper deals with the problem of existence of a convergent ``strong'' normal form in the neighbourhood of an equilibrium, for a finite dimensional system of differential equations with analytic and 
time-dependent non-linear term. The problem can be solved either under some \emph{non-resonance} hypotheses on the spectrum of the linear part or if the non-linear term is assumed to be (slowly) decaying in time. This paper ``completes'' a pioneering work of Pustil'nikov in which, despite under weaker non-resonance hypotheses, the nonlinearity is required to be asymptotically autonomous. The result is obtained as a consequence of the existence of a strong normal form for a suitable class of real-analytic Hamiltonians with non-autonomous perturbations. 
 
\end{abstract}

\section{Preliminaries and main result}
\subsection{Introduction}
The study of the dynamics of a system of ordinary differential equations (ODEs) in a neighbourhood of an equilibrium, boasts nowadays a rich and well established theory. Its foundation goes even back to the late XIX century to the contribution of Poincar\'{e} \cite{poi79} and Lyapunov \cite{lyap92}. Given an analytic vector field, the possibility to write the motions of the associated system in the vicinity of an equilibrium as a convergent power series, is deeply related to some \emph{non-resonance} conditions on the eigenvalues of the linear part.\\   
The results have been afterwards extended in the studies of Siegel started in \cite{sieg42}. The problem of the reducibility of a given system to a linear form via an analytic transformation, it is shown to be solvable in \cite{siegel52} for a full measure set of eigenvalues . \\
In the case of Hamiltonian structure, investigated later in \cite{siegel54} , the problem can be naturally interpreted in terms of the existence of a (convergent) canonical transformation of variables, casting a Hamiltonian of the form ``quadratic'' $+$ ``perturbation'' into a suitable\footnote{I.e. such that the corresponding canonical equations are integrable.} normal form, in some neighbourhood of the examined equilibrium. Based on this approach, the paper \cite{giorLyap} provides a generalisation of the results by Lyapunov, removing the hypothesis of purely imaginary eigenvalues. \\
In any case, we remark that, as a common feature of this class of problems, without any assumption on the eigenvalues, the program of casting the Hamiltonian at hand into a normal form, at least in general, fails. In fact, it is immediate to recognize how the linear combinations of eigenvalues occurring in the normalization scheme could produce some ``small divisor effects''. Knowingly, this phenomenon can either obstruct the formal resolvability of the homological equations produced during the normalization or jeopardize the convergence of the series. \\
We recall that, for instance, the described problem of well-posedness of the homological equation is overcome by Moser in \cite{moser56}, in the case of ``one and a half\footnote{With periodic time dependence.}'' degrees of freedom Hamiltonian $H(p,q,t)$ close to a hyperbolic equilibrium located at $p=q=0$. The strategy consists of keeping terms of the form $(pq)^k$, $k \geq 2$, in the normal form. In this way the canonical equations are still integrable ($x:=pq$ is a prime integral) but this allows to avoid the division by zero in the homological equation which would have been carried by those terms. This analysis plays a fundamental role in the context of instability phenomena in Hamiltonian systems with several degrees of freedom (Arnold's diffusion), in order to describe the flow in the neighbourhood of partially hyperbolic tori of a priori unstable systems, see \cite{chga94}. \\
The pioneering work by Pustil'nikov \cite{pust}, aims to extend the results of the paper \cite{siegel52}, by introducing a time dependence in the non-linear part of the vector field (not necessarily Hamiltonian). As it is natural, the choice of a suitable class of time-dependent perturbations and its treatment is a further difficulty to the phenomenon of the  ``resonances''. In \cite{pust}, under the non-resonance condition already assumed in \cite{siegel52} for the autonomous case, it is required that the perturbation is asymptotic to a time-independent, analytic function. However, no restrictions are imposed on the ``type'' of the time dependence, more specifically, it has to be neither periodic nor quasi-periodic. This case is also known as \emph{aperiodic} time dependence.\\
After \cite{pust}, the interest in a general dependence on time has been renewed in \cite{giozen} then followed by \cite{boun13}, \cite{fw14a} and subsequent papers. Basically, all of them deal with the Hamiltonian case (see \cite{fw15a} for the case of Poisson systems). The paper \cite{fw15} extends the above described result by Moser to the case of a perturbation aperiodically dependent on time. \\
As a matter of fact, the Hamiltonian structure is not a real obstruction for the use of the tools apt to treat the Hamiltonian case. In fact, given a system of ODEs, it can be always interpreted as (``a half'' of the) canonical equations of a suitable Hamiltonian system, of larger dimension, see e.g. \cite{berdvar}. The strategy of this paper is to derive the integrability of the system of ODEs at hand, see \ff{eq:sys}, as a particular case of the existence of a normal form for a  real-analytic Hamiltonian with aperiodic perturbation, see \ff{eq:ham}, by using the tools introduced in \cite{fw15} for the one degrees of freedom case.\\
The possibility to cast the Hamiltonian \ff{eq:ham} into a normal form is shown to be possible in the two cases described in Theorem \ref{thm}. In the second case, we deal with perturbations linear in the $y$ variables, in the presence of some non-resonance assumption on the eigenvalues. This case is directly related to the Hamiltonian formulation of a system of ODEs (due to the linearity in $y$). It is immediate to notice that, with respect to \cite[(0.3)]{pust}, the condition \ff{eq:nonres} on the eigenvalues is clearly  more restrictive. Nevertheless, the hypothesis of asymptotic time-independence assumed in \cite{pust} is weakened to the simple boundedness.\\   
On the other hand, the first case, has a more general character: if the perturbation decays\footnote{The exponential decay, see \ff{eq:decay}, is chosen for simplicity of discussion. The only necessary assumption is the summability in $t$ of the perturbing function over the non-negative real semi-axis, see \cite{fw15c}.} in time, either the described assumption on the form of $f$ or on the eigenvalues turn out to be unnecessary. Basically, the presence of resonance phenomena is no longer an obstruction for the existence of the normal form, see also  \cite{fw15c}.\\
The paper, based on the \emph{Lie series} formalism developed by A. Giorgilli et al., can be regarded, at the same time, as a non-autonomous version of \cite{giorLyap}.

\subsection{Setting}
Let us consider the following Hamiltonian 
\beq{eq:ham}
H(x,y,\eta,t)=h(x,y,\eta)+ f(x,y,t), \qquad h(x,y,\eta):=\eta + \sum_{l=1}^n \lambda_l x_l y_l \mx{,}
\eeq 
where $(x,y,\eta) \in \ml{D}:=[-r,r]^n \x [-r,r]^n \x \RR$, with $n \geq 1$ and $r>0$,  $\lambda_l \in \CC$ and $t \in \RR^+:=[0,+\infty)$. The assumptions on $f$ will be discussed below. The system (\ref{eq:ham}) is nothing but the ``autonomous equivalent'' of $\ml{H}(x,y,t)=\sum_{l=1}^n \lambda_l x_l y_l  + f(y,x,t)$, once $\eta$ has been defined as the conjugate variable to $t$.\\
The standard use of the analytic tools requires the complexification of the domain $\ml{D}$ as follows. Given $R \in (0,1/2]$ set $\ml{D}_R:=\ml{Q}_R \times \ml{S}_{R} $, where
$$
\ml{Q}_R:=\{(x,y) \in \CC^{2n} : |x|,|y| \leq  R \},\qquad \ml{S}_R:=\{\eta \in \CC : |\Im \eta| \leq  R \} \mx{,}
$$  
It will be required that, for all $t \in \RR^+$, $f$ belongs to the space of real-analytic functions on $\accentset{\circ}{\ml{Q}}_R$ and continuous on the boundary, which we denote with $\mathfrak{C}(\ml{Q}_R)$. In such a way $ H \in \mathfrak{C}(\ml{D}_R) $.\\In particular, the space of all the $G \in \mathfrak{C}(\ml{Q}_R) $ is endowed with the 
\emph{Taylor norm}
\beq{eq:taylor}
\norm{G(x,y,t)}{R}:=\sum_{\alpha,\beta \in \NN^n} |g_{\alpha,\beta}(t)| R^{|\alpha+\beta|} \mx{,} 
\eeq
where $G(x,y,t)=:\sum_{\alpha,\beta \in \NN^n} g_{\alpha,\beta}(t) x^{\alpha} y^{\beta}$ and\footnote{It is understood that $x^{\alpha}y^{\beta}:=x_1^{\alpha_1}\cdot\ldots\cdot x_n^{\alpha_n} \cdot y_1^{\beta_1}\cdot\ldots\cdot y_n^{\beta_n}$.} $|\alpha|:=\sum_{l=1}^n \alpha_l$. We recall the standard result for which, if $G \in \mathfrak{C}(\ml{Q}_R)$ for all $t \in \RR^+$, then  $|g_{\alpha,\beta}(t)| \leq \snorm{G}{R}R^{-|\alpha+\beta|}$, where $|G|_R:=\sup_{(x,y) \in \ml{Q}_R}|G|$. In particular, $\norm{G}{R'}<+\infty$ for all $R'<R$.\\
Throughout this paper we shall deal with perturbations satisfying the following conditions:
\begin{enumerate}
\item $f$ is ``at least'' quadratic in $x$ and ``at least'' linear in $y$: a property that we will denote with ($QxLy$), i.e. $f_{\alpha,\beta}(t) = 0$ for all $t \in \RR^+$ and for all $(\alpha,\beta) \in \NN^{2n} \setminus  \Gamma$, where $\Gamma:=\{(\alpha,\beta) \in \NN^{2 n}:|\alpha| \geq 2, \, |\beta| \geq 1 \}$,
\item there exist $M_f  \in [1,+\infty)$ and $a \in [0,1)$ such that\footnote{The interval $a \in [0,1)$ is a compact way to denote either the time decay $a \in (0,1)$ or the boundedness $a=0$. As in our previous paper we recall that we are interested in the case of small $a$ (slow decay) and the upper bound $a=1$ is set for simplicity. On the other hand, it is easy to realise that the case $a \geq 1$ is straightforward.}, for all $(x,y,t) \in \ml{Q}_R \times \RR^+$, 
\beq{eq:decay}
\norm{f(x,y,t)}{R} \leq M_f e^{-a t} \mx{.}
\eeq
\end{enumerate}
\subsection{Main result}
In the described setting, the main result can be stated as follows
\begin{satz}\label{thm} 
Suppose that one of the following conditions are satisfied: 
\begin{description}
\item{I. Time decay:} $a>0$.
\item{II. Linearity in $y$ $+$ non-resonance:} $a=0$ and the perturbation is linear in $y$, denoted by ($Ly$), i.e. of the form $f(x,y,t)=y \cdot g(x,t)$. In addition, the vector $\Lambda:=(\lambda_1,\ldots,\lambda_n)$, satisfies the \emph{non-resonance condition}
\beq{eq:nonres}
\max_{l=1,\ldots,n} \left(|\Re \ml{U}(\alpha,e_l,\Lambda)|^{-1}  \right)\leq \gamma |\alpha|^{\tau}, \qquad \forall  \alpha \in \NN^n \mx{,}
\eeq  
where $\ml{U}(\alpha,\beta,\Lambda):=(\alpha-\beta) \cdot \Lambda$, for some $\gamma>0$ and $\tau \geq n$. $e_l$ stands for the $l-$th vector of the canonical basis of $\RR^n$.
\end{description}
Then it is possible to determine $R_*,R_0$ with $0<R_*<R_0 \leq R^{16}$ and a family of canonical transformations $(x,y,\eta)=\ml{M}(x^{(\infty)},y^{(\infty)},\eta^{(\infty)})$, $\ml{M}:\ml{D}_{R_*} \rw \ml{D}_{R_0}$, analytic on $\ml{D}_{R_*}$ for all $t \in \RR^+$, casting the Hamiltonian (\ref{eq:ham}) into the \emph{strong normal form} 
\beq{eq:normalform}
H^{(\infty)}(x^{(\infty)},y^{(\infty)},\eta^{(\infty)})=h(x^{(\infty)}, y^{(\infty)} , \eta^{(\infty)})  \mx{.}
\eeq
\end{satz}
\begin{rem}
It is immediate to recognize the similarity between (\ref{eq:nonres}) and the standard Diophantine condition. Clearly, all the vectors $\Lambda$ whose real part is a Diophantine vector, satisfy condition (\ref{eq:nonres}), no matter what the imaginary part is. Hence the set of vectors satisfying (\ref{eq:nonres}) is, \emph{a fortiori}, a full-measure set.\\
As anticipated in the introduction, we stress that condition (\ref{eq:nonres}) is stronger than the non-resonance condition imposed in \cite{pust} and it is not satisfied in the case of purely imaginary $\Lambda$.
\end{rem}
\begin{rem}\label{remtwo}
As usually done in the \emph{Lie series method}, see e.g. \cite{gio03}, the transformation $\ml{M}$ will be constructed as the limit (defined, at the moment, only at a formal level)  
\beq{eq:composition}
\ml{M}:=\lim_{j \rw \infty} \ml{M}^{(j)} \circ \ml{M}^{(j-1)} \circ \ldots \circ \ml{M}^{(0)} \mx{,}
\eeq
where $\ml{M}^{(j)}:=\exp(\lie{j}) \equiv \id + \sum_{s \geq 1} (s!)^{-1} \lie{j}^s$ and $\lie{j}:=\{\cdot,\chi^{(j)}\}$. The \emph{generating sequence} $\{\chi^{(j)}\}_{j \in \NN}$, where $\chi^{(j)}=\chi^{(j)}(x,y,t)$, see \cite{giozen}, is meant to be determined. \\
We will show (see the proof of Lemma \ref{lem}) that in the case of a perturbation which is ($Ly$), it is possible to show that $\chi^{(j)}(x,y,t)$ is ($Ly$) as well, for all $j \in \NN$. In such a case, it is easy to check by induction that $x\jj= \ml{M}^{(j)} x^{(j+1)}$ \emph{does not} depend on the variable $y$, for all $j$. Hence the composition $x \equiv x^{(0)}=\ml{M} x\ift =:\ml{M}_x (x\ift,t)$ does not depend on $y\ift$ i.e. is an analytic map $\ml{M}_x:  \tilde{\ml{Q}}_{R_*} \rw \tilde{\ml{Q}}_{R_0}$ parametrised by $t$, where $\tilde{\ml{Q}}_{R}:=\{x \in \CC^n : |x|\leq R\}$.
This will play a key role in the next section. 
\end{rem}

\subsection{The corollary}
Let us consider the following non-linear system
\beq{eq:sys}
\dot{v}=Av+g(v,t) \mx{,}
\eeq
where $v \in \RR^n$, $A$ is a $n \x n$ matrix with real entries and  the function $g$ is such that $\pl_{v}^{\nu}g(0,t) \equiv 0$ for all $\nu \in \NN^n$ such that $|\nu| \leq 1$ i.e. $g$ is at least quadratic in $v$. We restrict ourselves to the class of diagonalizable $A$ with non-purely imaginary eigenvalues $\lambda_l$. In the obvious system of coordinates denoted with $x$, the system (\ref{eq:sys}) easily reads as 
\beq{eq:sysdiag}
\dot{x}_l=\lambda_l x_l + \tilde{g}_l(x,t), \qquad l=1,\ldots,n\mx{.}
\eeq
In this framework one can state the next
\begin{cor}
Suppose that  $f(x,y,t):=y \cdot \tilde{g}(x,t)$ and $\Lambda$ is such that the conditions described in II of Theorem \ref{thm} are satisfied. Then the system (\ref{eq:sysdiag}) is integrable in a suitable neighbourhood of the origin. \\ 
The same result holds, in particular, without any non-resonance condition on $\Lambda$, provided that $\tilde{g}(x,t)$ is such that (\ref{eq:decay}) is satisfied with $a>0$.
\end{cor}
\proof 
The key remark, see e.g. \cite{berdvar}, is that (\ref{eq:sysdiag}) can be interpreted as a set of canonical equations of the Hamiltonian system with Hamiltonian $
\ml{K}:=\eta + \sum_{l=1}^n y_l(\Lambda_l x_l + \tilde{g}_l(x,t))$, i.e. (\ref{eq:ham}) with $f(x,y,t)$ defined in the statement. Hence, by Theorem \ref{thm}, there exists a suitable neighbourhood of the origin endowed with a set of coordinates $(x\ift,y\ift,\eta\ift)$, such that $\ml{K}$ is cast into the (integrable) strong form $\ml{K}\ift=\eta\ift+ \sum_{l=1}^n \lambda_l y_l\ift x_l\ift$. Furthermore, as noticed in Remark \ref{remtwo}, $\ml{M}_x$ is an analytic map between $x$ and $x\ift$. Hence $x(t)=\ml{M}_x (x\ift(0)\exp(\ml{A} t),t)$, with $\ml{A}:=\diag(\lambda_1,\ldots,\lambda_n)$, gives the explicit solution of (\ref{eq:sysdiag}). 
\endproof

%%%%%%%%%%%%%%%%%%%%%%%%%%%%%%%%%%%%%%%%%%%%%%%%%%%%%%%%%%%%%%%%%%%%%%%%%%%%%%%%%%%%
%% SECTION THE SECOND %%%%%%%%%%%%%%%%%%%%%%%%%%%%%%%%%%%%%%%%%%%%%%%%%%%%%%%%%%%%%%
%%%%%%%%%%%%%%%%%%%%%%%%%%%%%%%%%%%%%%%%%%%%%%%%%%%%%%%%%%%%%%%%%%%%%%%%%%%%

\section{Some preliminary results}
\subsection{Two elementary inequalities}
\begin{prop}\label{prop:estimates}
For all $\ml{R} \leq e^{-4}$ and all $\delta \leq 1/2$ the following inequalities hold
\beq{eq:in}
\sum_{\substack{\nu \in \NN^m \\ |\nu|\geq N}} \ml{R}^{|\nu|} \leq 2 m e^{3m-3} \ml{R}^{\frac{3N}{4}},\qquad
\sum_{\nu \in \NN^m} |\nu|^{\mu} (1-\delta)^{|\nu|} \leq \ml{C}(m,\mu) \delta^{-m-\mu-1} \mx{,}
\eeq 
where $m \geq 2$, $\mu \geq 0$ and $\ml{C}(m,\mu):=e^{4m+\mu-1} (m+\mu)^{(m+\mu)}/(m-1)!$.
\end{prop}
\proof See Appendix.
\endproof

\subsection{A result on the homological equation}
\begin{prop}\label{propone}
Consider the following equation
\beq{eq:homol}
\lie{j} h + f\jj=0 \mx{,} 
\eeq
where $h$ has been defined in (\ref{eq:ham}) and $f\jj=f\jj(x,y,t)=\sum_{(\alpha,\beta) \in \Gamma} f_{\alpha,\beta}\jj (t) x^{\alpha} y^{\beta}$  satisfies $\norm{f\jj}{\tilde{R}} \leq M_j \exp(-a t)$ for some $a \in [0,1)$. The following statements hold for all $\delta \in (0,1/2]$: 
\begin{enumerate}
\item If $a>0$, there exists $C_1=C_1(n,\Lambda)>0$ such that  
\beq{eq:estimone}
\norm{\chi\jj}{(1-\delta)\tilde{R}},\norm{\pl_t \chi\jj}{(1-\delta)\tilde{R}} \leq C_1 M_j a^{-1} \delta^{-2(n+1)} \mx{.}
\eeq
\item If $a=0$, $f\jj$ is of the form $f\jj=y \cdot g\jj (x,t)$ and $\Lambda$ satisfies (\ref{eq:nonres}), there exists $C_2=C_2(n,\Lambda,\tau,\gamma)>0$ such that
\beq{eq:estimtwo}
\norm{\chi\jj}{(1-\delta)\tilde{R}},\norm{\pl_t \chi\jj}{(1-\delta)\tilde{R}} \leq C_2 M_j \delta^{-(n+\tau+2)} \mx{.}
\eeq
\end{enumerate}
\end{prop}
\proof
First of all note that $\lie{j} h=\pl_t \chi \jj+ \sum_{l=1}^n \lambda_l (x_l \pl_{x_l}-y_l \pl_{y_l})\chi\jj  $. By expanding the generating function as $\chi\jj (x,y,t)=\sum_{(\alpha,\beta) \in \NN^{2n}} c_{\alpha,\beta}\jj (t) x^{\alpha} y^{\beta}$, equation (\ref{eq:homol}) reads, in terms of Taylor coefficients, as 
\beq{eq:homolocomp}
\dot{c}_{\alpha,\beta}\jj (t) + \ml{U}(\alpha,\beta,\Lambda) c_{\alpha,\beta}\jj =f_{\alpha,\beta}\jj (t) \mx{.} 
\eeq
The solution of (\ref{eq:homolocomp}) is easily written, for all $(\alpha,\beta) \in \Gamma$, as
\beq{eq:solhom}
c_{\alpha,\beta}\jj (t) = e^{-\ml{U}(\alpha,\beta,\Lambda)t}\left[c_{\alpha,\beta}\jj (0) + \int_0^t e^{\ml{U}(\alpha,\beta,\Lambda)s} f_{\alpha,\beta}\jj (s) ds  \right] \mx{,}
\eeq 
while trivially $c_{\alpha,\beta}\jj(t) \equiv 0$ for all $(\alpha,\beta) \in \NN^{2n} \setminus \Gamma$.\\ 
Now denote $\ml{U}_R+i \ml{U}_I:=\ml{U}(\alpha,\beta,\Lambda)$ with $\ml{U}_{I,R} \in \RR$ and recall that, by hypothesis, $|f_{\alpha,\beta}\jj (t)| \leq M_j \tilde{R}^{-|\alpha+\beta|} e^{-a t}$.\\
\textbf{Case} $a>0$. For all $(\alpha,\beta) \in \Gamma$ such that $\ml{U}_R \geq 0$ we choose $c_{\alpha,\beta}\jj (0)=0$ then we have 
\[|c_{\alpha,\beta}\jj| \leq e^{-\ml{U}_R t} \int_0^t e^{\ml{U}_R s} |f_{\alpha,\beta}\jj (s)|ds \leq M_j \tilde{R}^{-|\alpha+\beta|} \int_0^t e^{-as} ds \leq M_j \tilde{R}^{-|\alpha+\beta|} a^{-1} \mx{.}\] 
Otherwise, for those $\alpha$ and $\beta$ such that $\ml{U}_R<0$, redefine $\ml{U}_R:=-\ml{U}_R$ with $\ml{U}_R>0$ and choose $c_{\alpha,\beta}\jj (0):=-\int_{\RR^+} \exp(\ml{U}(\alpha,\beta,\Lambda)s) f_{\alpha,\beta}\jj(s) ds$. Note that $|c_{\alpha,\beta}\jj (0)| < + \infty$. In this case we have $|c_{\alpha,\beta}\jj| \leq \exp(\ml{U}_R t) \int_t^{\infty} \exp(-\ml{U}_R s) |f_{\alpha,\beta}\jj (s)|ds \leq M_j \tilde{R}^{-|\alpha+\beta|} a^{-1}$. Hence $|c_{\alpha,\beta}\jj|\leq M_j \tilde{R}^{-|\alpha+\beta|} a^{-1}$ for all $(\alpha,\beta) \in \Gamma$. By recalling (\ref{eq:taylor}) one gets $\norm{\chi\jj}{(1-\delta)\tilde{R}} \leq M_j a^{-1} \sum_{(\alpha,\beta) \in \NN^{2n}} (1-\delta)^{|\alpha+\beta|}$. The use of the second of \ff{eq:in} with $\nu:=(\alpha,\beta)$, yields the first part of (\ref{eq:estimone}) with $C_1$ set for the moment to $\hat{C}_1:=\ml{C}(2n,0)$.\\
Directly from (\ref{eq:homolocomp}) we get $|\dot{c}_{\alpha,\beta}\jj| \leq |\alpha+\beta||\Lambda||c_{\alpha,\beta}\jj|+|f_{\alpha,\beta}\jj| \leq  a^{-1} M_j (1+|\Lambda|)  |\alpha+\beta|\tilde{R}^{-|\alpha+\beta|}$. By \ff{eq:in} with $\mu=1$ we get the second of part of \ff{eq:estimone}. The constant is chosen as $C_1:=(1+ |\Lambda|) \ml{C}(2n,1)  > \hat{C}_1$. \\
\textbf{Case} $a=0$. In such case, the homological equation reads as 
\beq{eq:homolocomptwo}
\dot{c}_{\alpha,l}\jj (t) + \ml{U}(\alpha,e_l,\Lambda) c_{\alpha,l}\jj =f_{\alpha,l}\jj (t) \mx{,} 
\eeq
where $f_{\alpha,l}\jj:=f_{\alpha,\beta}\jj|_{\beta=e_l}$
(the same notation for $c_{\alpha,l}\jj$), for all $\alpha \in \NN^n$ such that $|\alpha| \geq 2$ and for all $l=1,\ldots,n$. By hypothesis (\ref{eq:nonres}), $\ml{U}_R \neq 0$. Similarly to the case $a>0$, if $\ml{U}_R >0$ we set $c_{\alpha,l}\jj (0) =0$, otherwhise, $c_{\alpha,l}\jj (0):=-\int_{\RR^+} \exp(\ml{U}(\alpha,e_l,\Lambda)s) f_{\alpha,l}\jj(s) ds$. Proceeding as before, one obtains, by using (\ref{eq:nonres}),
\[
|c_{\alpha,l}\jj (t)| \leq M_j \ml{U}_R^{-1} \tilde{R}^{-|\alpha|-1} \leq \gamma M_j |\alpha|^{\tau} \tilde{R}^{-|\alpha|-1} \mx{.}
\]
This implies $\norm{\chi\jj}{(1-\delta)\tilde{R}} \leq n \gamma M_j \sum_{\alpha \in \NN^n} |\alpha|^{\tau} (1-\delta)^{|\alpha|}$ which is, by \ff{eq:in}, the first part of \ff{eq:estimtwo} with $\hat{C}_2=n \gamma \ml{C}(n,\tau)$. On the other hand, from the homological equation, we get $|\dot{c}_{\alpha,l}\jj(t)| \leq M_j |\alpha|^{\tau+1} (1+\gamma|\Lambda|) \tilde{R}^{-|\alpha|-1}$. Similarly, the latter yields the second part of \ff{eq:estimtwo} with $C_2:=\max\{ n (1+\gamma|\Lambda|) \ml{C}(n,\tau+1), \hat{C}_2\}$. 
\endproof

\subsection{A bound on the Lie operator}
\begin{prop}\label{proptwo}
Let $F,G$ be two functions such that $\norm{F}{(1-\tilde{d})\tilde{R}},\norm{G}{(1-\tilde{d})\tilde{R}}< +\infty$ for some $\tilde{d} \in (0,1/4]$ and $\tilde{R}>0$. Then for all $ s \in \NN$ the following bound holds
\beq{eq:spoisson}
\norm{\ml{L}_{G}^s F}{(1-2 \tilde{d})\tilde{R}} \leq e^{-2} s! [e^2 (\tilde{R} \tilde{d} )^{-2} \norm{G}{(1-\tilde{d})\tilde{R}}]^s \norm{F}{(1-\tilde{d})\tilde{R}} \mx{.}
\eeq
\end{prop}
\proof
Straightforward from \cite[Sec 3.2]{giorLyap} and \cite[Lemma 4.2]{gio03}.
\endproof

%%%%%%%%%%%%%%%%%%%%%%%%%%%%%%%%%%%%%%%%%%%%%%%%%%%%%%%%%%%%%%%%%%%%%%%%%%%%%%%%%%%%
%% SECTION THE THIRD %%%%%%%%%%%%%%%%%%%%%%%%%%%%%%%%%%%%%%%%%%%%%%%%%%%%%%%%%%%%%%%
%%%%%%%%%%%%%%%%%%%%%%%%%%%%%%%%%%%%%%%%%%%%%%%%%%%%%%%%%%%%%%%%%%%%%%%%%%%%

\section{Proof of the main result: convergence of the normal form}
\subsection{Preparation of the domains}
Taking into account the domain restriction imposed by Proposition \ref{proptwo}, the canonical transformations will be constructed of the form $\ml{M}_j:\ml{D}_{R_{j+1}} \rw \ml{D}_{R_j} \ni (x\jj,y\jj,\eta\jj)$ (understood $(x^{(0)},y^{(0)},\eta^{(0)}) \equiv (x,y,\eta)$), where $\{\ml{D}_{R_j}\}_{j \in \NN}$ is a suitable sequence of nested domains. We will also provide another sequence $\{\epsilon_j\}$ which will be used to control the size of the remainder.     
\begin{lem}\label{propthree}
Let us consider the following sequences 
\beq{eq:rec}
\epsilon_{j+1}=K a^{-1} d_j^{-\sigma} \epsilon_j^{2},\qquad 
R_{j+1}:=(1-2 d_j) R_j \mx{,}
\eeq
with $\epsilon_j,R_j < 1$, $d_j \leq 1/4$ and where $\epsilon_0, R_0,a,K,\sigma>0$ are given.  If
\beq{eq:convcond}
\epsilon_0 \leq \epsilon_a:=a (2 \pi)^{-\sigma} K^{-1} \mx{,}
\eeq
then it is possible to construct $\{d_j\}_{j \in \NN}$ in such a way $R_j \geq R_*:=R_0/2$ and $\epsilon_j \rw 0$ monotonically as $j \rw \infty$.  
\end{lem}
\begin{rem}
The property $R_*>0$ is crucial, as $R_*$ is the lower bound for the analyticity radius of the normalised Hamiltonian.
\end{rem}
\proof Straightforward from \cite[Lemma 4.4]{fw15c}. We recall that a suitable choice is $\epsilon_j=\epsilon_0(j+1)^{-\sigma}$, then, by \ff{eq:rec}, $d_j=(\epsilon_0  K a^{-1} )^{(1/\sigma)}(j+2)^2/(j+1)^4$. From the latter, one has 
\beq{eq:series}
\sum_{j \geq 0} d_j \leq 1/6 \mx{,} 
\eeq
provided that condition \ff{eq:convcond} is satisfied.
\endproof
\subsection{Iterative lemma}
Let us define for all $j \geq 0$, $H^{(j+1)}:=\ml{M}_j H\jj$ with $H^{(0)}:=H$.  
\begin{lem}\label{lem}
Under the same hypotheses of Theorem \ref{thm} and under the condition \ff{eq:convcond} it is possible to find a $R_0$ and a sequence $\{\chi\jj\}_{j \in \NN}$ such that $H^{(j)}(x,y,\eta,t)=h(x,y,\eta)+f\jj (x,y,t)$ with $f\jj$ ($QxLy$) and such that $\norm{f\jj}{R_j} \leq \epsilon_j e^{-at}$ for all $j$, where $\epsilon_j,R_j$ are given by \ff{eq:rec}.
\end{lem} 
The stated result exploits the possibility to remove the perturbation with the normalization algorithm obtaining, in this way, the desired normal form \ff{eq:normalform}. The interpretation of $\epsilon_j$ as a bound for the remainder is clearly related to the well known feature of the \emph{quadratic method}.
\proof
 By induction. If $j=0$, the statement is clearly true by hypothesis, by setting $f^{(0)}:=f$, either in the case $I$ or in the case $II$. We are supposing here that $\epsilon_0$ is small enough in order to satisfy \ff{eq:convcond}. This will be achieved later by a suitable choice of $R_0$.\\
Let us suppose the statement to be valid for $j$. In this way we get
\[
H^{(j+1)} \equiv \exp(\lie{j}) H\jj  =   h+ f\jj + \lie{j} h + \sum_{s \geq 1} (s!)^{-1} \lie{j}^s f\jj + \sum_{s \geq 2} (s!)^{-1} \lie{j}^s h \mx{.} 
\]
We shall determine $\chi\jj$ in such a way \ff{eq:homol} is satisfied so that, by setting 
\beq{eq:fjpo}
f^{(j+1)}:=\sum_{s \geq 1} \frac{1}{s!} \lie{j}^s f\jj + \sum_{s \geq 2} \frac{1}{s!} \lie{j}^s h \Heq{\ref{eq:homol}}{=} \sum_{s \geq 1} \frac{s}{(s+1)!} \lie{j}^s f\jj \mx{,}
\eeq 
one has $H^{(j+1)}=h+f^{(j+1)}$.\\
It is immediate from \ff{eq:homolocomp} that $\chi\jj$ has the same null Taylor coefficients as $f\jj$. Hence if $f\jj$ is ($QxLy$) then $\chi\jj$ is also. It is easy to check by induction that this implies that $\lie{j}^s f\jj$ is ($QxLy$) for all $s$, then $f^{(j+1)}$ is ($QxLy$). Similarly, equation \ff{eq:homolocomptwo} implies that if $f\jj$ is ($Ly$) then $\chi\jj$ is also. This implies that $\lie{j}^s f\jj $ is ($Ly$) for all $s$, hence $f^{(j+1)}$ is ($Ly$). This completes the formal part. In particular, by induction, $f\jj $ is ($Ly$) for all $j$, as claimed in Remark \ref{remtwo}. \\ 
Let us now discuss the quantitative estimate on $f\jj$ in the case $a>0$. By Propositions \ref{propone}, \ref{proptwo} and the inductive hypothesis, one gets 
\beq{eq:lief}
\norm{\lie{j}^s f\jj}{(1-2 d_j)R_j} \leq s! \Theta^s  \epsilon_j e^{-at},\qquad
\Theta:=\frac{e^2  C_1 }{ a R_*^2 d_j^{2n+4}} \epsilon_j \mx{.}
\eeq
Setting $K:=2 n e^2 C_1 R_*^{-2}$ and $\sigma:=2n + 5$, we have that 
\beq{eq:}
2 n \Theta=(K \epsilon_j a^{-1} d_j^{-\sigma})d_j  \leq d_j  \mx{,}
\eeq
as $\epsilon_{j+1}/\epsilon_j<1$ by Lemma \ref{propthree}. Hence, $\Theta < 1/2$ and the series defined in \ff{eq:fjpo} is convergent, furthermore
\beq{eq:lastlief}
e^{at} \norm{f^{(j+1)}}{R_{j+1}} \leq \epsilon_j \sum_{s \geq 1} \Theta^s  \leq  2 n \Theta \epsilon_j \Heq{\ref{eq:}}{\leq}  K a^{-1} d_j^{-\sigma} \epsilon_j^2 \Heq{\ref{eq:rec}}{=}  \epsilon_{j+1} \mx{,}
\eeq
which completes the inductive step. The condition \ff{eq:convcond} in this case reads as 
\beq{eq:condone}
\epsilon_0 \leq a R_0^2 (2\pi)^{-\sigma} (8 n e^2 C_1)^{-1} \mx{.}
\eeq
On the other hand, from the analyticity of $f$, we get $|f_{\alpha,\beta}(t)| \leq M_f R^{-|\alpha+\beta|} \leq M_f R_0^{-|\alpha+\beta|/16}$, as $R_0 \leq R^{16}$ by hypothesis. By using the first of \ff{eq:in} we get $\norm{f}{R_0} \leq M_f \sum_{(\alpha,\beta) \in \NN^{2n}} R_0^{(15/16)|\alpha+\beta|} \leq 2n e^{(2n-1)} M_f R_0^{135/64}=:\epsilon_0$. Replacing the latter in \ff{eq:condone}, the condition on $R_0$ described in the statement of Theorem \ref{thm} is meant to be completed with the following one 
\beq{eq:rzeroa}
R_0  \leq [a/(16 (2 \pi)^{\sigma} e^{2n+1} n^2 C_1 M_f) ]^{64/7} \mx{.}
\eeq
The case $a=0$ is analogous: it is sufficient to replace $C_1$ with $C_2$, remove the term $e^{\pm a t}$ from the statement, \ff{eq:lief} and \ff{eq:lastlief}, then replace $a$ with $1$ from \ff{eq:lief} to \ff{eq:condone}, where now $\sigma=n+\tau+5$. The only substantial difference consists in the sum obtained from \ff{eq:in}, which is slightly improved, since $f$ linear in $y$. We have in this case $\norm{f}{R_0} \leq n^2 e^{n-1} M_f R_0^{75/32}=:\epsilon_0$ leading to 
\beq{eq:rlesszeroa}
R_0  \leq [8 (2 \pi)^{\sigma} e^{n+1} n^3 C_2 M_f]^{-32/11}\mx{.}
\eeq
\endproof

\subsection{Bounds on the coordinate transformation}
\begin{lem}
The transformation of coordinates defined by the limit \ff{eq:composition} satisfies
\beq{eq:trasf}
|x^{(\infty)}-x|,|y^{(\infty)}-y|,|\eta^{(\infty)}-\eta| \leq R_0/6 \mx{,}
\eeq 
in particular, it defines an analytic map $\ml{M}:\ml{D}_{R_*} \rw \ml{D}_{R_0}$ and $H^{(\infty)}:=\ml{M} H$ is an analytic function on $\ml{D}_{R_*}$.
\end{lem}
\proof
We will discuss the case $a>0$. The case $a=0$ is straightforward simply replacing $C_1$ with $C_2$, $a$ with $1$ and changing the value of $\sigma$, where necessary.\\
Let us start from the variable $x$. Note that, by Proposition \ref{proptwo}, one has $\norm{\lie{j}^s x_l^{(j+1)}}{(1-2 d_j)R_j} \leq s! \Theta^s R_0$ for all $l=1,\ldots,n$. Hence we have, by \ff{eq:}
\[
|x^{(j+1)}-x^{(j)}| \leq n \max_{l=1,\ldots,n} \sum_{s \geq 1} \frac{1}{s!}
\norm{\lie{j}^s x_l^{(j+1)}}{(1-2 d_j)R_j} \leq 2 n R_0 \Theta \leq R_0 d_j \mx{.}
\]
In this way $ |x^{(\infty)}-x| \leq \sum_{j \geq 0} |x^{(j+1)}-x^{(j)}| $ converges by \ff{eq:series}. The procedure for $y$ is analogous. \\
As for the third of \ff{eq:trasf}, it is necessary to observe that $\lie{j} \eta = - \pl_t \chi^{(j)}$. Hence, by \ff{eq:spoisson} and the second of \ff{eq:estimone}, one has $
\norm{\lie{j}^s \eta}{(1-2 R_j)} \leq e^{-2}s! \Theta^{s-1} (R_*^2 e^{-2}\Theta) \leq s! \Theta^s R_0$, hence $|\eta^{(j+1)}-\eta^{(j)}| \leq 2 n R_0 \Theta \leq R_0 d_j $.\\
The bounds \ff{eq:trasf} ensure that points in $\ml{D}_{R_*}$ are mapped within $\ml{D}_{R_0}$ where $R_*=R_0/2$. Furthermore, the absolute convergence of the above described series, ensured by \ff{eq:series}, guarantees the uniform convergence in every compact subset of $\ml{D}_{R_*}$ and the analyticity of $\ml{M}$, and then of $H^{(\infty)}$, follows from the theorem of Weierstra\ss, see e.g. \cite{dett}.
\endproof

\subsection*{Appendix. Proof of Proposition \ref{prop:estimates}}
First of all, recall $\sum_{|\nu| \geq N} |\nu|^{\mu} \ml{R}^{|\nu|} = 
\sum_{l \geq N} 
\binom{l+m-1}{m-1} l^{\mu} \ml{R}^l$. 
Now note that $\log \prod_{j=1}^{m-1}(l+j) \leq \int_1^m \log(l+x)dx =1-m+\log[(m+l)^{(m+l)}(1+l)^{-(1+l)}]$ hence $(m-1)!\binom{l+m-1}{m-1}= \prod_{j=1}^{m-1}(l+j) \leq e^{m-1} (m+l)^{(m+l)}(1+l)^{-(1+l)} \leq e^{2m-2} (m+l)^{(m+\mu)} $. This yields 
\beq{eq:intin}
\sum_{|\nu| \geq N} |\nu|^{\mu} \ml{R}^{|\nu|} \leq [e^{2m-2}/(m-1)!]
\sum_{l \geq N} (m+l)^{(m+\mu)} \ml{R}^{l} \mx{.}
\eeq
On the other hand, the function $h(x):=(m+x)^{\kappa}\ml{R}^{x/4}$ has a maximum in $x=0$ (in the non-negative semi-axis) if $\ml{R} \leq \exp(-4 \kappa /m)$ and in $x^*:=-m-4 \kappa/\log\ml{R} $ otherwise. Hence, from (\ref{eq:intin}) with $\mu=0$ we have $\sum_{|\nu| \geq N} \ml{R}^{|\nu|} \leq [(m-1)!]^{-1} m^m e^{2m-2}
\sum_{l \geq N} \ml{R}^{(3/4)l}$ which gives the first of \ff{eq:in} by using the inequality $m^m \leq e^{m-1}m!$ and recalling $\ml{R} \leq e^{-4}$.\\
Now set $\ml{R}=1-\delta$. By hypothesis $\ml{R}>e^{-4}$, hence $(m+l)^{(m+\mu)}(1-\delta)^{l/4} \leq (1-\delta)^{-m/2}(-2 (m+\mu)/\log (1-\delta))^{(m+\mu)}$. By substituting the latter in \ff{eq:intin} with $N=0$, then using the inequalities $-\log(1-\delta) \geq \delta$ and $[1-(1-\delta)^{3/4}] \geq \delta/2$ as $\delta \leq 1/2$, the second of  \ff{eq:in} easily follows.
\subsection*{Acknowledgements}
The first author is grateful to Prof. Dario Bambusi for remarkable discussions on this problem. 

\bibliographystyle{alpha}
\bibliography{Pst.bib}

\newcommand{\sort}[1]{}
\begin{thebibliography}{FW15b}

\bibitem[Ber09]{berdvar}
V.~Berdichevsky.
\newblock {\em Variational Principles of Continuum Mechanics: I. Fundamentals}.
\newblock Interaction of Mechanics and Mathematics. Springer Berlin Heidelberg,
  2009.

\bibitem[Bou13]{boun13}
A.~Bounemoura.
\newblock Effective stability for slow time-dependent near-integrable
  hamiltonians and application.
\newblock {\em C. R. Math. Acad. Sci. Paris}, 351(17-18):673--676, 2013.

\bibitem[CG94]{chga94}
L.~Chierchia and G.~Gallavotti.
\newblock Drift and diffusion in phase space.
\newblock {\em Ann. Inst. H. Poincar\'e Phys. Th\'eor.}, 60(1):144, 1994.

\bibitem[Det65]{dett}
J.W. Dettman.
\newblock {\em Applied Complex Variables}.
\newblock A Series of advanced mathematics texts. Dover Publications, 1965.

\bibitem[FW14]{fw14a}
A.~Fortunati and S.~Wiggins.
\newblock Normal form and {N}ekhoroshev stability for nearly integrable
  {H}amiltonian systems with unconditionally slow aperiodic time dependence.
\newblock {\em Regul. Chaotic Dyn.}, 19(3):363--373, {\sort{a}} 2014.

\bibitem[FW15a]{fw15a}
A.~Fortunati and S.~Wiggins.
\newblock {A Kolmogorov theorem for nearly-integrable Poisson systems with
  asymptotically decaying time-dependent perturbation.}
\newblock {\em Regul. Chaotic Dyn.}, 20(4):476--485, 2015.

\bibitem[FW15b]{fw15}
A.~Fortunati and S.~Wiggins.
\newblock Normal forms à la {M}oser for aperiodically time-dependent
  {H}amiltonians in the vicinity of a hyperbolic equilibrium.
\newblock Accepted for the publication on Discr. Cont. Dyn. Sys. -S., 2015.

\bibitem[FW15c]{fw15c}
A.~Fortunati and S.~Wiggins.
\newblock Negligibility of small divisor effects in the normal form theory for
  nearly-integrable hamiltonians with decaying non-autonomous perturbations.
\newblock \url{http://arxiv.org/abs/1509.02119}, {\sort{c}} 2015.

\bibitem[Gio]{giorLyap}
A.~Giorgilli.
\newblock On a {T}heorem of {L}yapounov.
\newblock {\em Rendiconti dell'Istituto Lombardo Accademia di Scienze e
  Lettere, Classe di Scienze Matematiche e Naturali}, In Press.

\bibitem[Gio03]{gio03}
A.~Giorgilli.
\newblock Exponential stability of {H}amiltonian systems.
\newblock In {\em Dynamical systems. {P}art {I}}, Pubbl. Cent. Ric. Mat. Ennio
  Giorgi, pages 87--198. Scuola Norm. Sup., Pisa, 2003.

\bibitem[GZ92]{giozen}
A.~Giorgilli and E.~Zehnder.
\newblock Exponential stability for time dependent potentials.
\newblock {\em Z. Angew. Math. Phys.}, 43(5):827--855, 1992.

\bibitem[Lya92]{lyap92}
A.M. Lyapunov.
\newblock {\em General Problem of the Stability Of Motion}.
\newblock Control Theory and Applications Series. Taylor \& Francis, 1992.

\bibitem[Mos56]{moser56}
J.~Moser.
\newblock The analytic invariants of an area-preserving mapping near a
  hyperbolic fixed point.
\newblock {\em Comm. Pure Appl. Math.}, 9:673--692, 1956.

\bibitem[Poi79]{poi79}
H.~Poincar{\'e}.
\newblock {\em Sur les propri{\'e}t{\'e}s des fonctions d{\'e}finies par les
  {\'e}quations aux diff{\'e}rences partielles}.
\newblock Paris. Universit{\'e}. Facult{\'e} des sciences. Theses.
  Gauthier-Villars, 1879.

\bibitem[Pus74]{pust}
L.~D. Pustyl'nikov.
\newblock Stable and oscillating motions in nonautonomous dynamical systems.
  {A} generalization of {C}. {L}. {S}iegel's theorem to the nonautonomous case.
\newblock {\em Mat. Sb. (N.S.)}, 94(136):407--429, 495, 1974.

\bibitem[Sie42]{sieg42}
C.~L. Siegel.
\newblock Iteration of analytic functions.
\newblock {\em Ann. of Math. (2)}, 43:607--612, 1942.

\bibitem[Sie52]{siegel52}
C.~L. Siegel.
\newblock \"{U}ber die {N}ormalform analytischer {D}ifferentialgleichungen in
  der {N}\"ahe einer {G}leichgewichtsl\"osung.
\newblock {\em Nachr. Akad. Wiss. G\"ottingen. Math.-Phys. Kl.
  Math.-Phys.-Chem. Abt.}, 1952:21--30, 1952.

\bibitem[Sie54]{siegel54}
C.~L. Siegel.
\newblock \"{U}ber die {E}xistenz einer {N}ormalform analytischer
  {H}amiltonscher {D}ifferentialgleichungen in der {N}\"ahe einer
  {G}leichgewichtsl\"osung.
\newblock {\em Math. Ann.}, 128:144--170, 1954.

\end{thebibliography}

\end{document}